\documentclass{amsart}
\usepackage{verbatim,amsmath,amssymb}

\newcommand{\fh}{\hat{f}}

\newcommand{\geh}{\mathfrak{g}}
\newcommand{\heh}{\mathfrak{h}}

\newcommand{\gc}{\geh_0}
\newcommand{\Pc}{P_0}

\newcommand{\Qc}{Q_0}

\newcommand{\beh}{\mathfrak{b}}
\newcommand{\Ih}{I\backslash\{0\}}
\newcommand{\la}{\lambda}
\newcommand{\La}{\Lambda}
\newcommand{\om}{\omega}
\newcommand{\omt}{\widetilde{\omega}}
\newcommand{\ve}{\varepsilon}
\newcommand{\vp}{\varphi}
\newcommand{\al}{\alpha}

\newcommand{\inner}[2]{\langle #1\,,\,#2\rangle}
\newcommand{\wt}{\mathrm{wt}}

\newcommand{\tM}{\widetilde{M}}

\newcommand{\WX}{\widetilde{W}}

\newcommand{\Wc}{W_0}

\newcommand{\Z}{\mathbb{Z}}

\newcommand{\IS}{I^0}

\newcommand{\st}{S}
\newcommand{\Aut}{\mathrm{Aut}}

\newcommand{\ww}[1]{w_0^{#1}}

\newcommand{\vn}{\emptyset}

\newtheorem{theorem}{Theorem}
\newtheorem{proposition}[theorem]{Proposition}
\newtheorem{lemma}[theorem]{Lemma}
\newtheorem{corollary}[theorem]{Corollary}

\theoremstyle{definition}
\newtheorem{remark}{Remark}
\newtheorem{assumption}{Assumption}
\newtheorem{example}{Example}

\numberwithin{equation}{section}
\numberwithin{theorem}{section}

\begin{document}

\title{Demazure structure inside Kirillov--Reshetikhin crystals}

\author[G.~Fourier]{Ghislain Fourier}
\address{Mathematisches Institut der Universit\"at zu K\"oln,
Weyertal 86-90, 50931 K\"oln, Germany}
\email{gfourier@mi.uni-koeln.de}

\author[A.~Schilling]{Anne Schilling}
\address{Department of Mathematics, University of California, One Shields
Avenue, Davis, CA 95616-8633, U.S.A.}
\email{anne@math.ucdavis.edu}
\urladdr{http://www.math.ucdavis.edu/\~{}anne}

\author[M.~Shimozono]{Mark Shimozono}
\address{Department of Mathematics, Virginia Tech, Blacksburg, VA 24061-0123, U.S.A.}
\email{mshimo@math.vt.edu}
\urladdr{http://www.math.vt.edu/people/mshimo/}

\thanks{\textit{Date:} May 2006}
\thanks{AS was supported in part by NSF grant DMS-0200774 and DMS-0501101.}
\thanks{MS was supported in part by NSF grant DMS-0401012.}

\begin{abstract}
The conjecturally perfect Kirillov-Reshetikhin (KR) crystals are
known to be isomorphic as classical crystals to certain Demazure
subcrystals of crystal graphs of irreducible highest weight modules
over affine algebras. Under some assumptions we show that the
classical isomorphism from the Demazure crystal to the KR crystal,
sends zero arrows to zero arrows. This implies that the affine
crystal structure on these KR crystals is unique.
\end{abstract}

\maketitle

\section{Introduction}

The irreducible finite-dimensional modules over a quantized affine
algebra $U'_q(\geh)$ were classified by Chari and
Pressley~\cite{CP:1995,CP:1998} in terms of Drinfeld polynomials. We
are interested in the subfamily of such modules which possess a
global crystal basis. Kirillov--Reshetikhin (KR) modules are
finite-dimensional $U'_q(\geh)$-modules $W^{r,s}$ that were
introduced in \cite{HKOTT:2002,HKOTY:1999}. It is expected that each
KR module has a crystal basis $B^{r,s}$, and that every irreducible
finite-dimensional $U'_q(\geh)$-module with crystal basis, is a
tensor product of the crystal bases of KR modules.

The KR modules $W^{r,s}$ are indexed by a Dynkin node $r$ of the
classical subalgebra (that is, the distinguished simple Lie
subalgebra) $\gc$ of $\geh$ and a positive integer $s$. In general
the existence of $B^{r,s}$ remains an open question. For type
$A_n^{(1)}$ the crystal $B^{r,s}$ is known to
exist~\cite{KKMMNN:1992a} and its combinatorial structure has been
studied~\cite{S:2002}. In many cases, the crystals $B^{1,s}$ and $B^{r,1}$ for
nonexceptional types are also known to exist and their combinatorics has been
worked out in~\cite{KKM:1994,KKMMNN:1992a} and~\cite{JMO:2000,Ko:1999},
respectively.

Viewed as a $U_q(\gc)$-module by restriction, $W^{r,s}$ is generally
reducible; its decomposition into $U_q(\gc)$-irreducibles was
conjectured in \cite{HKOTT:2002,HKOTY:1999}. This was verified by
Chari~\cite{C:2001} for the nontwisted cases.

Kashiwara~\cite{Ka:2002} conjectured that as classical crystals,
many of the KR crystals (the ones conjectured to be perfect in
\cite{HKOTT:2002,HKOTY:1999}) are isomorphic to certain Demazure
subcrystals of affine highest weight crystals. Kashiwara's
conjecture was confirmed by Fourier and Littelmann~\cite{FL:2004} in
the untwisted cases and Naito and Sagaki~\cite{NS:2005} in the
twisted cases.

In this paper we prove that the classical isomorphism from the
Demazure crystals to KR crystals sends zero arrows to zero arrows
(see Theorem~\ref{thm:Dem KR}). It is not an affine crystal
isomorphism but becomes an isomorphism after tensoring with an
appropriate affine highest weight crystal. This recovers some of the
isomorphisms given by the Kyoto path model. We emphasize this is
accomplished without the assumption of perfectness of the KR
crystals. The automorphisms on the crystals that are used in the
definition of the ground state path in the Kyoto path model, come from
affine Dynkin diagram automorphisms which can be calculated using the
factorization of a translation element in the extended affine Weyl group in
our setting. For the proof of our results we require the assumptions of regularity of KR
crystals, the existence and uniqueness of a certain special element $u$
in a KR crystal, and the existence of automorphisms on KR crystals
coming from certain Dynkin automorphisms (see
Assumption~\ref{A:KR}). We show that under these assumptions, the KR
crystals admit a unique affine crystal structure (see
Corollary~\ref{cor:KR unique}), and we give an algorithm which shows
that twofold tensor products of KR crystals are connected (see
Corollary~\ref{cor:connect}). We expect that Assumption~ \ref{A:KR}
holds, that is, if the existence of the KR crystals were established
these hypotheses could be removed.

In Section~\ref{S:notation} we establish notation and review some
results about the extended affine Weyl group. The definition of
Demazure crystals and KR crystals is given in
Section~\ref{S:crystal}. Section~\ref{S:Dem KR} contains our main
result stated in Theorem~\ref{thm:Dem KR} showing that all zero
arrows of the Demazure crystal are present in the KR crystal. In
Section~\ref{S:hw} we provide explicit sequences of lowering
operators leading from the special element $u$ of a KR crystal to
all classical highest weight elements of the KR crystal. The
connectedness of tensor products of KR crystals and an application
regarding the algorithmic calculation of the combinatorial
$R$-matrix can be found in Section~\ref{S:connect}.

\subsection*{Acknowledgments} We like to thank Philip Sternberg
for many stimulating discussion regarding KR crystals of type
$D_n^{(1)}$.

\section{Notation and basics} \label{S:notation}

\subsection{Affine Kac-Moody algebras}
Let $\geh$ be an affine Kac-Moody algebra with Cartan subalgebra
$\heh$, Dynkin node set $I=\{0,1,\dotsc,n\}$, Cartan matrix
$A=(a_{ij})_{i,j\in I}$, realized by the set of linearly independent
simple roots $\{\alpha_i\mid i\in I\} \subset \heh^*$ and simple
coroots $\{\al_i^\vee\mid i\in I\} \subset \heh$, such that
$\inner{\al_i^\vee}{\al_j}=a_{ij}$ \cite{Kac:1990}. Let $d\in\heh$
be the scaling element, which is any element such that
$\inner{d}{\al_i}=0$ for $i\in \Ih$ and $\inner{d}{\al_0}=1$. Let
$(a_i\mid i\in I)$ be the unique tuple of relatively prime positive integers
that give a linear dependence relation among the columns of $A$, and
let $(a_i^\vee\mid i\in I)$ be the tuple for the rows of $A$. Let
$\delta=\sum_{i\in I} a_i \al_i$ be the null root,
$\theta=\sum_{i\in \Ih} a_i\al_i$, and $c=\sum_{i\in I} a_i^\vee
\al_i^\vee$ the canonical central element. We
have $\inner{d}{\delta}=a_0$. Let $\{\La_i\mid i\in I\}\subset
\heh^*$ be the fundamental weights, which, together with
$\delta/a_0$, are defined to the dual basis to the basis
$\{\al_i^\vee\mid i\in I\} \cup \{d\}$ of $\heh$. In particular
$\inner{\al_i^\vee}{\La_j}=\delta_{ij}$. Let $P=\bigoplus_{i\in I}
\Z \La_i \oplus \Z(\delta/a_0)\subset \heh^*$ be the weight lattice,
$P^+=\bigoplus_{i\in I} \Z_{\ge0} \La_i\oplus
\Z(\delta/a_0)=\{\la\in P\mid \inner{\al_i^\vee}{\la}\ge0 \text{ for
all $i\in I$}\}$ the set of dominant weights and $Q=\bigoplus_{i\in
I} \Z \al_i\subset \heh^*$ the root lattice. The level of a weight
$\la\in P$ is defined by $\inner{c}{\la}$. Let $W$ be the affine
Weyl group, generated by the simple reflections $\{s_i\mid i\in
I\}$. $W$ acts on $P$ by $s_i\la = \la - \inner{\al_i^\vee}{\la}\al_i$.

Let $(\cdot\mid \cdot)$ be the nondegenerate $W$-invariant
symmetric form on $\heh^*$; it is defined by
$(\al_i\mid \al_j) = a_i^\vee a_i^{-1} a_{ij}$ for $i,j\in I$,
$(\al_i\mid \La_0)=0$ for $i\in\Ih$, $(\al_0\mid \La_0)=a_0^{-1}$,
and $(\La_0\mid \La_0)=0$. One may check that
\cite[(6.4.1)]{Kac:1990}
\begin{align} \label{E:thetalen}
 (\theta\mid \theta) = 2a_0 = \begin{cases}
   4 &\text{for $A_{2n}^{(2)}$} \\
   2 &\text{otherwise.}
 \end{cases}
\end{align}
The pairing $(\cdot\mid \cdot)$ induces an isomorphism
$\nu:\heh\rightarrow \heh^*$ given by $\inner{\nu(h)}{h'}=(h\mid
h')$ for all $h,h'\in \heh$. So $\nu(\al_i^\vee)=a_i
(a_i^\vee)^{-1} \al_i$ for $i\in I$, $\nu(d)=a_0\La_0$, and
$\nu(c)=\delta$. Define $\theta^\vee\in\heh$ by
$\nu(\theta^\vee)=2\theta/(\theta\mid\theta)=\theta/a_0$.

Let $\gc\subset \geh$ be the simple Lie subalgebra whose Dynkin node
set is $\Ih$, with Weyl group $\Wc\subset W$, root lattice $\Qc$,
weight lattice $\Pc$, and fundamental weights $\{\om_i\mid
i\in\Ih\}\subset \Pc$.

Let $P'=P/\Z(\delta/a_0)$. The natural projection $P'\rightarrow
\Pc$ has a section $\Pc\rightarrow P'$ defined by $\om_i\mapsto
\La_i-a_i^\vee\La_0$ for $i\in \Ih$. The image of this section is
the set of elements in $P'$ of level zero.

\subsection{Dynkin automorphisms} \label{SS:Dynkinaut}
Let $X$ denote the affine Dynkin diagram and $\Aut(X)$ denote the
group of automorphisms of $X$. By definition an element of $\Aut(X)$
is a permutation of the Dynkin node set $I$ which preserves the kind
of bonds between nodes. Observe that
\begin{equation} \label{E:auta}
\begin{aligned}
  a_{\tau(i)} &= a_i \\
  a_{\tau(i)}^\vee &= a_i^\vee
\end{aligned} \qquad\qquad\text{for all $i\in I$ and $\tau\in\Aut(X)$.}
\end{equation}
There is an action of $\Aut(X)$ on $P$ given by
\begin{align*}
  \sigma(\La_i) &= \La_{\sigma(i)} \qquad\text{for $i\in I$} \\
  \sigma(\delta) &= \delta
\end{align*}
for $\sigma\in\Aut(X)$. By \eqref{E:auta} this action restricts to
an action of $\Aut(X)$ on $\Pc$ called the level zero action.

\subsection{Translations} For $\alpha\in\Pc$,
define the element $t_\alpha\in\Aut(P)$ by \cite[(6.5.2)]{Kac:1990}
\begin{align} \label{E:trans}
  t_\alpha(\la) =
  \la+\inner{c}{\la}\alpha-\bigl((\la\mid \alpha)+\dfrac{1}{2}(\alpha\mid\alpha)
  \inner{c}{\la}\bigr)\delta.
\end{align}
The map $\alpha\mapsto t_\alpha$ defines an injective group
homomorphism $\Pc\rightarrow\Aut(P)$ whose image shall be denoted
$T(\Pc)$. For any $w\in \Wc$,
\begin{align} \label{E:W0trans}
  w t_\alpha w^{-1} = t_{w(\alpha)}.
\end{align}
Therefore $\Wc \ltimes T(\Pc)$ acts on $P$. There is an induced
action of $\Wc\ltimes T(\Pc)$ on $P'$ that preserves the level of a
weight. For every $m\in\Z$ there is an action of $\Wc \ltimes
T(\Pc)$ on $\Pc$ called the level $m$ action, given by $w
*_m \mu = w (m\La_0+\mu)-m\La_0$ for $\mu\in \Pc$. Under the level
one action, the element $t_\alpha$ is precisely translation by
$\alpha$.

\subsection{Extended affine Weyl group}
For each $i \in \Ih$, define $c_i = \max(1,a_i/a_i^{\vee})$; these
constants were introduced in \cite{HKOTT:2002}. Using the Kac
indexing of the affine Dynkin diagrams \cite[Table Fin, Aff1 and Aff2]{Kac:1990},
we have $c_i=1$ except for $c_i=2$ for $\geh=B_n^{(1)}$ and $i=n$,
$\geh=C_n^{(1)}$ and $1\le i\le n-1$, $\geh=F_4^{(1)}$ and $i=3,4$,
and $c_2=3$ for $\geh=G_2^{(1)}$. Consider the sublattices
of $\Pc$ given by
\begin{align*}
  M &= \bigoplus_{i\in\Ih} \Z c_i\al_i = \Z \Wc \cdot \theta/a_0 \\
  \tM &= \bigoplus_{i\in\Ih} \Z c_i\om_i.
\end{align*}
It is easy to check that $M\subset \tM$ and that the action of $\Wc$
on $\Pc$ restricts to actions on $M$ and $\tM$. Let $T(\tM)$ (resp.
$T(M)$) be the subgroup of $T(\Pc)$ generated by $t_\la$ for
$\la\in\tM$ (resp. $\la\in M$).

%\remind{This table can be deleted.}
%\begin{align*}
%\begin{array}{|c|c|c|c|c|} \hline
%\geh & \gc & \theta/a_0 & M & \tM \\ \hline \hline %
%A_n^{(1)} & A_n & \eps_1-\eps_{n+1} & Q(A_n) & P(A_n) \\ %
%B_n^{(1)} & B_n & \eps_1+\eps_2 & Q(C_n) & P(C_n) \\ %
%C_n^{(1)} & C_n & 2\eps_1 & 2 Q(B_n) & 2 P(B_n) \\ %
%D_n^{(1)} & D_n & \eps_1+\eps_2 & Q(D_n) & P(D_n) \\ %
%A_{2n}^{(2)} & C_n & \eps_1 & P(C_n) & P(C_n) \\ %
%A_{2n-1}^{(2)} & C_n & \eps_1+\eps_2 & Q(C_n) & P(C_n) \\ %
%D_{n+1}^{(2)} & B_n & \eps_1 & Q(B_n) & P(B_n) \\ \hline
%\end{array}
%\end{align*}

There is an isomorphism \cite[Prop. 6.5]{Kac:1990}
\begin{align} \label{E:Waffiso}
  W \cong \Wc \ltimes T(M)
\end{align}
as subgroups of $\Aut(P)$. Under this isomorphism we have
\begin{align} \label{E:s0}
  s_0 = t_{\theta/a_0} s_\theta.
\end{align}
Define the extended affine Weyl group to be the subgroup of
$\Aut(P)$ given by
\begin{align} \label{E:WXdef}
 \WX = \Wc \ltimes T(\tM).
\end{align}
When $\geh$ is of untwisted type, $M\cong Q^\vee$, $\tM\cong
P^\vee$, $c_i\om_i = \nu(\om_i^\vee)$, and
$c_i\alpha_i=\nu(\alpha_i^\vee)$ for $i\in\Ih$.

Let $C\subset P \otimes_{\Z} \mathbb{R}$ be the fundamental chamber,
the set of elements $\la$ such that $\inner{\al^\vee_i}{\la}\ge0$
for all $i\in I$. Define the subgroup $\Sigma\subset\WX$ to be the
set of elements that send $C$ into itself.

It follows from \eqref{E:W0trans} and \eqref{E:Waffiso} that $W$ is
a normal subgroup of $\WX$. Thus $\Sigma$ acts on $W$ by
conjugation. Since the Weyl chambers adjacent to $C$ are precisely
those of the form $s_i(C)$ for $i\in I$, the element $\tau\in\Sigma$
induces a permutation (also denoted $\tau$) of the set $I$ given by
\begin{align}\label{E:autoW}
  \tau s_i \tau^{-1} = s_{\tau(i)} \qquad\text{for $i\in I$.}
\end{align}
Since the braid relations in $W$ are preserved, $\Sigma$ is a
subgroup of $\Aut(X)$.

\subsection{Special automorphisms} We identify the subgroup $\Sigma$
explicitly. Say that an affine Dynkin node $i\in I$ is
\textit{special} if there is an automorphism $\tau\in \Aut(X)$ of
the affine Dynkin diagram such that $\tau(i)=0$. In the untwisted
case, $i$ is special if and only if $\om_i^\vee$ is a minuscule
coweight. Let $\IS\subset I$ denote the set of special vertices.
Explicitly, using the Kac labeling~\cite{Kac:1990}:
\begin{align*}
  \IS = \begin{cases}
    \{0,1,\dotsc,n\} & \text{for $A^{(1)}_n$} \\
    \{0,1\} &\text{for $B^{(1)}_n$, $A^{(2)}_{2n-1}$} \\
    \{0,n\} &\text{for $C^{(1)}_n$, $D^{(2)}_{n+1}$} \\
    \{0,1,n-1,n\} &\text{for $D^{(1)}_n$} \\
    \{0,1,5\}&\text{for $E^{(1)}_6$} \\
    \{0,6\} &\text{for $E^{(1)}_7$} \\
    \{0\} &\text{otherwise.}
  \end{cases}
\end{align*}

\begin{proposition} \label{P:auto} For each $i\in\IS$ there is a
unique element $\tau_i\in\Sigma$ such that $\tau_i(i)=0$. Moreover
$\Sigma=\{\tau_i\mid i\in \IS\}$.
\end{proposition}

We call $\tau_i$ the special automorphism associated with $i\in
\IS$.

Note that every Dynkin automorphism is determined by its action on
$\IS$. We describe the special automorphisms explicitly. $\tau_0$ is
the identity automorphism. If $\geh$ is of untwisted affine type and
$i\in \IS$ then for all $j\in \IS$, $\tau_i(j)=k\in \IS$ where
$-\om_i+\om_j\cong \om_k \mod \Qc$ and $\om_0=0$ by convention. For
$\geh$ of twisted type the only nonidentity (special) automorphisms
are the elements of $\Aut(X)$ which on $\IS$ are given by
$\tau_1=(0,1)$ in type $A_{2n-1}^{(2)}$ and $\tau_n=(0,n)$ in type
$D^{(2)}_{n+1}$.

We now specify $\Sigma$ explicitly as a subgroup of permutations of
$\IS$. In all cases but $D^{(1)}_n$ and $n$ even, $\Sigma$ is a
cyclic group. This determines $\tau_i$ and $\Sigma$ completely
except for types $A^{(1)}_n$ and $D^{(1)}_n$. For $A^{(1)}_n$,
$\Sigma\cong \Z/(n+1)\Z$ where $\tau_i(j)=j-i\mod (n+1)$ for all
$i,j\in\IS$. For $D^{(1)}_n$ and $n$ odd, $\Sigma$ is cyclic with
$\tau_{n-1}=(0,n,1,n-1)$, $\tau_1=(0,1)(n-1,n)$ and
$\tau_n=(0,n-1,1,n)$ in cycle notation acting on $\IS$. For $n$
even, $\Sigma\cong \Z/2\Z \times \Z/2\Z$ with $\tau_1=(0,1)(n-1,n)$,
$\tau_{n-1}=(0,n-1)(1,n)$ and $\tau_n=(0,n)(1,n-1)$.

\begin{proposition} \label{P:WX} \
$\Sigma \cong \tM/M$ via $\tau_i\mapsto \om_i + M$ for $i\in\IS$
and
\begin{align} \label{E:extaffsigma}
  \WX \cong W \rtimes \Sigma.
\end{align}
as subgroups of $\Aut(\Pc)$.
\end{proposition}

If $i\in \IS$ then $c_i=1$ and we have
\begin{align} \label{E:specialWSigma}
\tau_i = \ww{\om_i} t_{-\om_i}
\end{align}
where, for $\la\in \Pc^+$,
\begin{align} \label{E:wwla}
\text{$\ww{\la}\in\Wc$ is the shortest element such that $\ww{\la}
\la$ is antidominant.}
\end{align}

\subsection{Dynkin automorphisms revisited}
Let $X_0$ be the Dynkin diagram for the classical subalgebra $\gc$
of $\geh$.

\begin{lemma} \label{L:F} There is a group homomorphism
\begin{equation} \label{E:autXX0}
\begin{split}
  \Aut(X) &\rightarrow\Aut(X_0) \\
  \sigma &\mapsto\sigma'
\end{split}
\end{equation}
where $\sigma'(i)=j$ if and only if $\sigma(\om_i)\in \Wc\om_j$.
\end{lemma}
\begin{proof} We first claim that there is a group action of
$\Aut(X)$ on $\Wc\backslash \Pc$ defined by $\sigma(\Wc \la) = \Wc
\sigma \la$ where $\Aut(X)$ acts on $\Pc$ via the level zero action.
 The level zero action of $s_0$ on $\Pc$ is the same as
that of $s_\theta\in \Wc$, by \eqref{E:s0} and \eqref{E:trans}. Thus
for the level zero action, $W \la = \Wc \la$ for $\la\in\Pc$. By
\eqref{E:autoW}, $\sigma \Wc \sigma^{-1} \subset W$ as it is
generated by $s_{\sigma(i)}$ for $i\in\Ih$. Thus we have
$\Wc\sigma\Wc \tau\la =\Wc (\sigma\Wc\sigma^{-1}) \sigma\tau\la =
\Wc \sigma\tau\la$. Therefore $\Aut(X)$ acts on $\Wc \backslash
\Pc$.

Next we show that this action restricts to an action on $F\subset
\Wc\backslash\Pc$ where $F$ is the set of $\Wc$-orbits of
fundamental weights $\om_i$ for $i\in\Ih$. Due to the above group
action we need only that $\sigma F\subset F$ for generators $\sigma$
of $\Aut(X)$. By \eqref{E:auta} we have $\sigma(\om_r)=
\om_{\sigma(r)} - a_r^\vee \om_{\sigma(0)}$ where we write
$\om_i=\La_i-a_i^\vee\La_0$ for all $i\in I$. Using this one may
straightforwardly check the lemma for each affine root system.
\end{proof}

$\Aut(X_0)$ is trivial except in the following cases, where the
homomorphism is described explicitly. The elements of $\Aut(X)$ and
$\Aut(X_0)$ are given by their action as permutations of $\IS$ and
$\IS\setminus\{0\}$ respectively.
\begin{enumerate}
\item $\Aut(A_n)$ is generated by the involution $i\mapsto n+1-i$
for $i\in\Ih$. In this case $\Aut(A_n^{(1)})$ is the dihedral group
$D_{2(n+1)}$. For $\sigma\in \Aut(A_n^{(1)})$, $\sigma'$ is the
nontrivial element in $\Aut(A_n)$ if and only if $\sigma$ reverses
orientation.
\item $\Aut(D_n)$ is generated by $(n-1,n)$ when $n>4$. In this case
$\Aut(D_n^{(1)})$ is generated by $(0,1)$, $(n-1,n)$ and
$(0,n)(1,n-1)$. All these map to the nontrivial element of
$\Aut(D_n)$ except in the case that $n$ is even, when $(0,n)(1,n-1)$
maps to the identity.
\item $\Aut(D_4)$ is the symmetric group on the three ``satellite"
vertices $\{1,3,4\}$. $\Aut(D_4^{(1)})$ is the symmetric group on
the vertices $\{0,1,3,4\}$ and is generated by $(0,i)$ for
$i\in\{1,3,4\}$. The generator $(0,i)$ is sent to the element
$(j,k)$ in $\Aut(D_4)$ where $\{0,i,j,k\}=\{0,1,3,4\}$ as sets.
\item $\Aut(E_6)$ is generated by
$(1,5)$. $\Aut(E_6^{(1)})$ is isomorphic to the $S_3$ that permutes
the special vertices $\{0,1,5\}$. Then each of the elements of order
two in $\Aut(E_6^{(1)})$ is sent to the nontrivial element of
$\Aut(E_6)$.
\end{enumerate}

\begin{remark}\label{R:autcl} In all cases, for all
$\tau\in\Sigma$, $\tau'$ is the identity in $\Aut(X_0)$. However for
$\sigma=(0,1)\in\Aut(D_n^{(1)})$ we have
$\sigma'=(n-1,n)\in\Aut(D_n)$.
\end{remark}

\section{Crystals} \label{S:crystal}

\subsection{Definition of crystals}
A $P$-weighted $I$-crystal is a set $B$, equipped with Kashiwara
operators $e_i,f_i:B\rightarrow B\sqcup \{\vn\}$, and weight
function $\wt:B\rightarrow P$ such that $e_i(f_i(b))=b$ if
$f_i(b)\not=\vn$, $f_i(e_i(b))=b$ if $e_i(b)\not=\vn$,
$\wt(f_i(b))=\wt(b)-\alpha_i$ if $f_i(b)\not=\vn$,
$\wt(e_i(b))=\wt(b)+\alpha_i$ if $e_i(b)\not=\vn$, and
$\inner{\al_i^\vee}{\wt(b)} = \vp_i(b)-\ve_i(b)$ where
$\vp_i(b)=\min\{ m\mid f_i^m(b)\not=\vn\}$ and $\ve_i(b)=\min\{m\mid
e_i^m(b)\not=\vn\}$ are assumed to be finite for all $b\in B$ and
$i\in I$. If $f_i(b)\not=\vn$ we draw an arrow colored $i$ from $b$
to $f_i(b)$. The connected components of the graph obtained by
removing all arrows from $B$ except the arrows colored $i$, are
called the $i$-strings of $B$. We write $\ve(b)=\sum_{i\in I}
\ve_i(b)\La_i$ and $\vp(b)=\sum_{i\in I} \vp_i(b)\La_i$.

An $I$-crystal $B$ is \textit{regular} if, for each subset $K\subset
I$ with $|K|=2$, each $K$-component of $B$ is isomorphic to the
crystal basis of an irreducible integrable highest weight
$U'_q(\geh_K)$-module where $\geh_K$ is the subalgebra of $\geh$
with simple roots $\alpha_i$ for $i\in K$.

The crystal reflection operator $\st_i:B\rightarrow B$ is defined by
the property that $\st_i(b)$ is the unique element in the $i$-string
of $b$ such that $\ve_i(\st_i(b))=\vp_i(b)$ or equivalently
$\vp_i(\st_i(b))=\ve_i(b)$. This defines an action of the Weyl group
$W$ on $B$ if $B$ is regular \cite{Ka:1994}.

If $B$ and $B'$ are $P$-weighted $I$-crystals, their tensor product
$B \otimes B'$ is a $P$-weighted $I$-crystal as follows (we use the
opposite of Kashiwara's convention). As a set $B\otimes B'$ is just
the Cartesian product $B \times B'$ where traditionally one writes
$b\otimes b'$ instead of $(b,b')$. The Kashiwara operators are given
by
\begin{align*}
f_i(b\otimes b') &= \begin{cases}
f_i(b)\otimes b'
& \text{if $\ve_i(b)\ge \vp_i(b')$} \\
b\otimes f_i(b') & \text{if $\ve_i(b)<\vp_i(b')$}
\end{cases} \\
e_i(b\otimes b') &=
\begin{cases} e_i(b) \otimes b' &
\text{if $\ve_i(b)>\vp_i(b')$} \\
b\otimes e_i(b') & \text{if $\ve_i(b)\le \vp_i(b')$.}
\end{cases}
\end{align*}

Given any $P$-weighted $I$-crystal $B$ and Dynkin automorphism
$\sigma$, there is a $P$-weighted $I$-crystal $B^\sigma$ whose
vertex set is written $\{b^\sigma \mid b\in B\}$ and whose edges are
given by $f_i(b)=b'$ in $B$ if and only if
$f_{\sigma(i)}(b^\sigma)=(b')^\sigma$. The weight function satisfies
$\wt(b^\sigma) =\sigma(\wt(b))$ where the second $\sigma$ is the
automorphism of $P$ defined by $\sigma$. A similar statement holds
for $\Pc$-weighted $I$-crystals, using the level zero action of
$\sigma$ on $\Pc$ defined in Subsection \ref{SS:Dynkinaut}.

Given any $P$-weighted $I$-crystal $B$, define the \textit{contragredient dual}
crystal $B^\vee = \{b^\vee\mid b\in B\}$ with $\wt(b^\vee)=-\wt(b)$
and $f_i(b)=b'$ if and only if $e_i(b^\vee)={b'}^\vee$.

\subsection{Branching} The following ideas have been applied
extensively (in \cite{KKMMNN:1992a} and \cite{StS:2004}, for
example) to identify the $0$-arrows in KR crystals. We shall use
them here for the same purpose.

Let $B$ be the crystal graph of a $U'_q(\geh)$-module and $K\subset
I$. A \textit{$K$-component} of $B$ is a connected component of the
graph obtained from $B$ by removing all $i$-edges for $i\not\in K$.
A \textit{$K$-highest weight vector} is an element $b\in B$ such
that $\ve_i(b)=0$ for all $i\in K$. Suppose $K$ is a proper subset
of $I$. Since the subalgebra of $\geh$ with simple roots
$\{\alpha_i\mid i\in K\}$ is semisimple, each $K$-component of $B$
has a unique $K$-highest weight vector. When $K=I\setminus\{0\}$ we
call the $K$-components and $K$-highest weight vectors
\textit{classical} components and highest weight vectors.

Suppose $\sigma$ is a Dynkin automorphism that fixes $K$ and induces
an automorphism (also denoted $\sigma$) on $B$ that sends $i$-arrows
to $\sigma(i)$-arrows for all $i\in I$. Then by definition $\sigma$
preserves $i$-arrows for all $i\in K$. There is a projection from
the classical weight lattice to that of the subalgebra with simple
roots $\alpha_i$ for $i\in K$; we refer to the latter as the
\textit{$K$-weight lattice}. In particular $\sigma$ permutes the
collection of $K$-components, sending $K$-highest weight vectors to
those with the same $K$-weight (that is, $\vp_i\circ\sigma=\vp_i$
for $i\in K$).

\subsection{Demazure modules and crystals}
Let $\geh$ be a symmetrizable Kac-Moody algebra and $U_q(\geh)$ its
quantized universal enveloping algebra. For a dominant weight $\La$
denote by $V(\La)$ the irreducible integrable highest weight
$U_q(\geh)$-module with highest weight $\La$. Write $B(\La)$ for its
crystal basis. Let $\beh$ be a Borel Lie subalgebra of $\geh$. For
$\mu\in W\cdot\La$ let $u_\mu$ be a generator of the line of weight
$\mu$ in $V(\La)$. Write $\mu=w\La$ where $w$ is shortest in its
coset $w W^\La$ and $W^\La=\{w\in W\mid w\La=\La\}$. When writing an
element $w\La\in W\cdot\La$ we shall always assume $w$ is of minimum
length. Define the Demazure module
$$
V_w(\La) := U_q(\beh)\cdot u_{w(\La)}.
$$
It is known that $V_w(\La)$ has a crystal base $B_w(\La)$
\cite{Ka:1993}; it is the full subgraph of $B(\La)$ whose vertex set
consists of the elements in $B(\La)$ that are reachable by raising
operators, from the unique element $u_{w\La}\in B(\La)$ of weight
$w\La$. We shall make use of the following result. By abuse of
notation let
\begin{align} \label{E:fw}
  f_w(b) = \{\,f_{i_N}^{m_N}\dotsm f_{i_1}^{m_1}(b)\mid
  m_k\in \Z_{\ge0}\}
\end{align}
where $w=s_{i_N}\dotsm s_{i_1}$ is any fixed reduced decomposition
of $w$. It is known~\cite{Ku:1987,L:1994,M:1988} that as sets,
\begin{align} \label{E:fwDem}
  B_w(\La) = f_w(u_{\La}).
\end{align}
For $\geh$ affine, let $w \in \widetilde{W}$. By
\eqref{E:extaffsigma} we may express it uniquely as $w = z\tau$
where $z \in W$ and $\tau \in \Sigma$. We define the Demazure module
to be
$$
V_{w}(\Lambda) := V_{z}(\tau(\Lambda)).
$$
Its crystal graph is denoted $B_w(\La)=B_z(\tau\La)$. For a
dominant $\la\in \tM$, let $\la^* = - w_0(\la)$, where $w_0$ is the
longest element in $W_0$. Define
$D(\la,s)=V_{t_{-\la^*}}(s\La_0)$ and by abuse of notation,
$D(\la,s)=B_{t_{-\la^*}}(s\La_0)$. For any $\sigma\in\Aut(X)$ let
$D^\sigma(\la,s) = B_{t_{-\sigma(\la)^*}}(s\La_{\sigma(0)})$; it
is obtained from $D(\la,s)$ by changing every $i$ arrow into a
$\sigma(i)$ arrow.

\subsection{KR crystals}
Kirillov--Reshetikhin (KR) modules $W^{r,s}$, labeled by $(r,s)\in \Ih\times \Z_{>0}$,
are finite-dimensional $U'_q(\geh)$-modules. See~\cite{HKOTT:2002} for the
precise definition. It is conjectured that $W^{r,s}$ has a global crystal 
basis $B^{r,s}$.

In~\cite{HKOTT:2002} a conjecture is given for the decomposition of
each Kirillov--Reshetikhin (KR) module $W^{r,c_r s}$ into its
$\geh_0$-components. Chari~\cite{C:2001} proved this conjecture for
the nonexceptional untwisted algebras and for the exceptional cases
for the nodes $r$ such that either $r\in \IS$ or $\om_r$ is the
highest root. Recently the $G_2$ case was treated in full
\cite{CM:2006}. In~\cite{FL:2004}, the $\geh_0$-structure of the
Demazure modules was calculated for the same
cases as in~\cite{C:2001}, and it was verified that the Demazure and
KR modules agree as $\geh_0$-modules. In addition, it was shown
in~\cite{FL:2005} that no matter what the precise $\geh_0$-structure
is, the Demazure and the KR modules agree as $\geh_0$-modules for
all untwisted algebras. Naito and Sagaki~\cite{NS:2005} proved the
conjectures of~\cite{HKOTT:2002} on the level of crystals for the
twisted cases under the assumption that the KR crystals for the
untwisted algebras exist. In unpublished work, Naito and Sagaki did
the same construction for the twisted cases on the Demazure modules.

\begin{remark} \label{R:KRDem} Assuming that $B^{r,c_rs}$ exists, the
Demazure crystal $D(c_r\om_r,s)$ and the KR crystal $B^{r,c_r s}$ have
the same classical crystal structure.
\end{remark}

In this paper we assume that the KR crystal $B^{r,c_rs}$ has the properties
of Assumption~\ref{A:KR}, which we expect to hold if the KR crystals exist.
In the next section we will see that with these assumptions the Demazure crystal
sits inside the KR crystal (see Theorem~\ref{thm:Dem KR}) and that the KR crystal
is unique (see Corollary~\ref{cor:KR unique}). For types $B_n^{(1)}$, $D_n^{(1)}$,
and $A_{2n-1}^{(2)}$ let $\sigma$ be the Dynkin automorphism exchanging the Dynkin
nodes $0$ and $1$ and fixing all others. For types $C_n^{(1)}$ and
$D_{n+1}^{(2)}$ let $\sigma$ be the Dynkin automorphism defined by
$i\mapsto n-i$ for all $i\in I$. We also write $\sigma$ for the
induced automorphism of $P$.

\begin{assumption} \label{A:KR} The KR crystal $B^{r,c_rs}$ has the following properties:
\begin{enumerate}
\item \label{A:regular} $B^{r,c_rs}$ is regular.
\item \label{A:u} There is a unique element $u\in B^{r,c_rs}$ such that
\begin{align*}
\ve(u)=s\La_0 \quad \text{and} \quad \vp(u) = s \La_{\tau(0)},
\end{align*}
where $t_{-c_r\om_r} = w \tau$ with $w\in W$ and $\tau\in \Sigma$.
\item \label{A:auto} For all types different from $A_{2n}^{(2)}$, $B^{r,c_r s}$ admits the
automorphism corresponding to $\sigma$ (also denoted $\sigma$) such that
\begin{equation}\label{E:sigmastring}
\ve\circ \sigma = \sigma\circ \ve \qquad \vp\circ\sigma
=\sigma\circ\vp.
\end{equation}
For type $A_{2n}^{(2)}$ we assume that $B^{r,c_rs}$ is given explicitly by the virtual
crystal construction in \cite{OSS:2003}.
\end{enumerate}
\end{assumption}

\section{Relation between Demazure and KR crystals} \label{S:Dem KR}

In this section we show that the Demazure crystal sits inside
the KR crystals in Theorem~\ref{thm:Dem KR} and, assuming their
existence, that the KR crystals are unique in Corollary~\ref{cor:KR unique}.

The main technique that we use in the proof is a decomposition of the
translation elements $t_{-c_r\om_r}$ that ends in a word for the subalgebra
associated to the nodes $\{0,1,\ldots,r-1\}$ of the Dynkin diagram in
analogy to the results of~\cite{FL:2004}.

\begin{proposition} \label{P:trans} Let $\geh$ be of nonexceptional
affine type, $r\in I\setminus\IS$ and $t_{-c_r\om_r} = w\tau$ for
$w\in W$ and $\tau\in\Sigma$. Then a reduced word for the minimum
length coset representative $w_2$ in $\Wc w$ is given by
\begin{equation} \label{E:wtilde}
w_2 =
  \begin{cases}
  \prod_{k=i}^1 s_0 (s_2s_3\dotsm s_{2k-1})(s_1s_2\dotsm s_{2k-2}) &
  \begin{matrix}
   \text{for $r=2i$ and} \\
   \text{$B_n^{(1)}$, $D_n^{(1)}$, $A_{2n-1}^{(2)}$}
  \end{matrix} \\
  \prod_{k=i}^1 s_0 (s_2s_3\dotsm s_{2k})(s_1s_2\dotsm s_{2k-1}) &
 \begin{matrix}
   \text{for $r=2i+1$ and} \\
   \text{$B_n^{(1)}$, $D_n^{(1)}$, $A_{2n-1}^{(2)}$}
  \end{matrix} \\
  \prod_{k=i}^1 s_0 (s_1 s_2\dotsm s_{k-1}) &
  \begin{matrix}
   \text{for $r=i$ and} \\
  \text{$C_n^{(1)}$, $A_{2n}^{(2)}$, $D_{n+1}^{(2)}$}
  \end{matrix}
  \end{cases}
\end{equation}
where the index $k$ decreases as the product is formed from left
to right.
\end{proposition}
\begin{proof} All nodes for $A_n^{(1)}$ are special so we may assume
$\geh$ is not of this type.

Applying the sequence of reflections in \eqref{E:wtilde} to
$\La_{\tau(0)}$, we see that each reflection $s_j$ changes the
weight by a positive multiple of $\alpha_j$, and the final weight is
$\La_0+c_r\om_r-i\delta$. It follows that \eqref{E:wtilde} yields a
reduced decomposition of some element $w_2\in W$.

Using \eqref{E:trans}, in all cases we have
\begin{align*}
  w\La_{\tau(0)} &= t_{-c_r \om_r} \tau^{-1} \La_{\tau(0)} = \La_0 - c_r\om_r -
  i\delta/a_0.
\end{align*}
Since $r\not\in \IS$ we have $w_0^{\om_r} \om_r =-\om_r$ where
$w_0^{\om_r}$ is defined in \eqref{E:wwla}. Moreover $w_0^{\om_r}$
is also the shortest element of $\Wc$ sending
$\La_0+c_r\om_r-i\delta/a_0$ to $\La_0-c_r\om_r-i\delta/a_0$. It
follows that $w = w_0^{\om_r} w_2$ is a length-additive
factorization and that $w_2$ is the minimum length coset
representative in $\Wc w$.
\end{proof}

\begin{remark} \label{R:w1w2} Let
$K=\{0,1,\dotsc,r-1\}\subset I$, $\geh_K\subset \geh$ the simple
subalgebra with Dynkin nodes $K$, $\{\omt_j\mid j\in K\}$ the
fundamental weights for $\geh_K$, and $W_K=\langle s_j\mid j\in
K\rangle\subset W$ the Weyl group of $\geh_K$. This given, we have
$w_2 = \ww{\omt_{\tau(0)}}$ where $\ww{\omt_j}\in W_K$ is defined
with respect to $\geh_K$.
\end{remark}

\begin{lemma} \label{L:simple} All of the weights of
$B^{r,c_rs}$ are in the convex hull of the $\Wc$-orbit $\Wc\cdot
c_rs\om_r$. Moreover for every $\mu\in \Wc\cdot c_rs\om_r$, there is
a unique element $u_\mu\in B(c_rs\om_r)\subset B^{r,c_rs}$ of the
extremal weight $\mu$.
\end{lemma}
\begin{proof} By \cite{FL:2004,NS:2005} the classical decomposition of
$D(c_r\om_r,s)$ agrees with that specified in \cite{HKOTT:2002}. In
every case the above condition holds.
\end{proof}

\begin{lemma} \label{lem:y} Let $\geh$ be of nonexceptional affine
type, $r\in I\setminus\IS$, $s\in\Z_{>0}$, $k<r$ where $B(c_r
s\om_k)$ occurs in $B^{r,c_rs}$, and $b=u_{c_r s\om_k}\in B(c_r
s\om_k)\subset B^{r,c_r s}$. Define
\begin{equation*}
y = \begin{cases} \st_2\dotsm \st_{k+1} \st_1\dotsm \st_k(b) &
\text{for $B_n^{(1)},D_n^{(1)},A_{2n-1}^{(2)}$,}\\
\st_1\dotsm\st_k(b) &\text{for $C_n^{(1)}, D_{n+1}^{(2)},
A_{2n}^{(2)}$.}
\end{cases}
\end{equation*}
Then
\begin{align} \label{E:f0s}
f_0^s(y) = \begin{cases} u_{c_r s\om_{k+2}} &\text{for $B_n^{(1)},
D_n^{(1)}, A_{2n-1}^{(2)}$,} \\
u_{c_r s\om_{k+1}} & \text{for $C_n^{(1)}, D_{n+1}^{(2)}, A_{2n}^{(2)}$.}
\end{cases}
\end{align}
\end{lemma}
\begin{proof} By definition the element $y$ is an extremal weight
vector within the classical crystal $B(c_r s\om_k)$. By weight
considerations one may check that
\begin{equation*}
y=\begin{cases} f_{2}^s  \cdots  f_{k}^s f_{k+1}^s f_{1}^s  f_{2}^s
\cdots
f_{k-1}^s f_{k}^s(b) & \text{for $B_n^{(1)},D_n^{(1)},A_{2n-1}^{(2)}$,}\\
f_1^{c_r s} f_2^{c_r s} \cdots f_k^{c_r s}(b) & \text{for
$C_n^{(1)}$, $D_{n+1}^{(2)}$, $A_{2n}^{(2)}$.}
\end{cases}
\end{equation*}
We claim that
\begin{equation*}
\begin{aligned}
&\ve(y)= s(\La_0+\La_2) && \vp(y)=s(\La_0+\La_{k+2})
&& \text{for $B_n^{(1)}$, $D_n^{(1)},A_{2n-1}^{(2)}$, $k>0$}\\
 &\ve(y)= s(\La_0+c_r\La_1) && \vp(y)=s(\La_0+c_r\La_{k+1}) &&
\text{for $C_n^{(1)}$, $A_{2n}^{(2)}$, $D_{n+1}^{(2)}$, $k>0$} \\
&\ve(y)=s\La_0 && \vp(y)=s\La_0&&\text{for $k=0$.}
\end{aligned}
\end{equation*}
By extremality and Lemma \ref{L:simple}, $y$ is in the indicated
position within its $i$-strings for $i\in\Ih$. It remains to show
that $\ve_0(y)=\vp_0(y)=s$ and \eqref{E:f0s} holds. In each case we
shall use Assumption \ref{A:KR} (\ref{A:auto}) either for the existence of a
crystal automorphism $\sigma$ on $B^{r,c_rs}$ or, in type
$A_{2n}^{(2)}$, for the virtual crystal construction of
$B^{r,c_rs}$.

We begin with type $D_n^{(1)}$. We have $c_r=1$ and
$\mu:=\wt(y)=(0^2,s^k,0^{n-k-2})$. Here we realize $\Pc\subset
((1/2)\Z)^n$ with $i$-th standard basis element $\epsilon_i$, with
$\om_i=(1^i,0^{n-i})$ for $1\le i\le n-2$ (we do not need the spin
weights) and $\al_i=\epsilon_i-\epsilon_{i+1}$ for $1\le i\le n-1$.
Let $b'=u_{s\om_{k+2}}\in B(s\om_{k+2})\subset B^{r,s}$. We have
$\vp_0(b')=0$, for otherwise $f_0(b')\in B^{r,s}$ has weight
contradicting Lemma \ref{L:simple}. Since
$\inner{\al_0^\vee}{\wt(b')}=2s$, we have $\ve_0(b')=2s$.

For type $D_n^{(1)}$, the automorphism $\sigma$ of $B^{r,c_rs}$
satisfies $e_0=\sigma\circ e_1 \circ\sigma$. Define
$z=e_1^s(\sigma(b'))$. It suffices to show that
\begin{align*}
  y =\sigma(z).
\end{align*}
Let $K=\{2,3,\dotsc,n\}\subset I$. The subalgebra of $\geh$ with
simple roots $\alpha_i$ for $i\in K$, is of type $D_{n-1}$. For this
reason we shall refer to $D_{n-1}$-components and $D_{n-1}$-highest
weight vectors instead of $K$-components and $K$-highest weight
vectors. Our proof rests on the following fact:
\begin{quote}
$B^{r,s}$ contains a unique element of weight $\mu$ that satisfies
$\ve_1=0$ and whose associated $D_{n-1}$-highest weight vector has
$D_n$-weight $\la:=(0,s^k,0^{n-k-1})$.
\end{quote}
For the classical components of $B^{r,s}$ that contain $D_{n-1}$-components
of weight $\la$, are precisely those of the form
$B((s-t)\om_k+t\om_{k+2})$ for $0\le t\le s$, and only for $t=0$
does the classical component contain an element of weight $\mu$ with
$\ve_1=0$ (and by extremality $B(s\om_k)$ contains a unique element
of weight $\mu$).

$y$ clearly satisfies the above property. It suffices to show that
$\sigma(z)$ does also.

$\sigma(b')$ is a $D_{n-1}$-highest weight vector with
$\wt(\sigma(b'))=(-s,s^{k+1},0^{n-k-2})$. So $\wt(z)=\mu$. By weight
considerations and Lemma \ref{L:simple},
$z'=\st_{k+1}\dotsm\st_2(z)$ is a $D_{n-1}$-highest weight vector of
weight $\la$. Therefore $\sigma(z)$ has weight $\sigma(\mu)=\mu$ and
has associated $D_{n-1}$-highest weight vector $\sigma(z')$, which
has weight $\sigma(\la)=\la$. Since the Dynkin nodes $0$ and $1$ are
nonadjacent we have $\ve_1(\sigma(z))=\ve_1(e_0^s(b'))=\ve_1(b')=0$.
Thus $\sigma(z)$ fulfills the above criteria and so must be equal to
$y$.

The proof is analogous for types $B_n^{(1)}$ and $A_{2n-1}^{(2)}$
using the same set $K$, which defines subalgebras of types $B_{n-1}$
and $C_{n-1}$ respectively.

For type $C_n^{(1)}$ we have $c_r=2$ for all $1\le r\le n-1$. Let
$K=\{1,2,\dotsc,n-1\}$; the associated subalgebra is of type
$A_{n-1}$. Here we realize $\Pc\cong\Z^n$ with $\om_i=(1^i,0^{n-i})$
for $1\le i\le n$ and $\al_i=\epsilon_i-\epsilon_{i+1}$ for $1\le
i\le n-1$ and $\al_n=2\epsilon_n$. Our argument uses the fact that
\begin{quote}
$B^{r,2s}$ contains a unique element of weight
$\mu:=(0,(2s)^k,0^{n-k-1})$ such that $\ve_n=0$ and whose associated
$A_{n-1}$-highest weight vector has $C_n$-weight $2s\om_k$.
\end{quote}
For the classical components in $B^{r,2s}$ that contain such an
$A_{n-1}$-component, are precisely those of the form
$B(2(s-t)\omega_k+2t\omega_{k+1})$ for $0\le t\le s$, and among
these, only for $t=0$ does the classical component contain an
element of weight $\mu$ for which $\ve_n=0$ (and by extremality
$B(2s\om_k)$ contains a unique element of weight $\mu$).

By construction $y$ satisfies this property. It suffices to show
that $\sigma(z)$ does also, where $z=e_n^s\circ\sigma(b')$ and
$b'=u_{2s\om_{k+1}}\in B(2s\om_{k+1})\subset B^{r,s}$.

We have $\vp_0(b')=0$ for otherwise $f_0(b')\in B^{r,2s}$ would have
weight contradicting Lemma \ref{L:simple}. Since
$\inner{\al_0^\vee}{\wt(b')}=2s$ we have $\ve_0(b')=2s$.

$\sigma(b')$ is an $A_{n-1}$-highest weight vector of weight
$\sigma(2s\om_{k+1})=(0^{n-k-1},(-2s)^{k+1})$. Therefore $z$ has
weight $(0^{n-k-1},(-2s)^k,0)$ and associated $A_{n-1}$-highest
weight vector $z'=\st_{n-k}\dotsm \st_{n-1}(z)$, which has weight
$(0^{n-k},(-2s)^k)$. It follows that $\sigma(z)$ has weight $\mu$
and its associated $A_{n-1}$-highest weight vector has weight
$2s\om_k$. Now $\ve_n(\sigma(z))=\ve_n(e_0^s(b'))=\ve_n(b')=0$ since
the Dynkin nodes $0$ and $n$ are nonadjacent. We have shown that
$\sigma(z)$ satisfies the above criteria and so must be equal to
$y$.

Type $D_{n+1}^{(2)}$ is similar to type $C_n^{(1)}$.

For type $A_{2n}^{(2)}$, the above kind of argument is not available
since $A_{2n}^{(2)}$ admits no nontrivial Dynkin automorphism.
Instead we apply virtual crystals. Under Assumption~\ref{A:KR}~(\ref{A:auto}), by
\cite{OSS:2003} the crystal $B^{r,s}$ is realized as the subset of
$V^{r,s}=B^{2n-r,s}_A \otimes B^{r,s}_A$ of type $A_{2n-1}^{(1)}$
generated from $u_{s\om_{2n-r}} \otimes u_{s\om_r}$ by the virtual
crystal operators $\fh_i=f_i f_{2n-i}$ for $1\le i\le n$ and
$\fh_0=f_0^2$ where $f_i$ are the crystal operators of the
$A_{2n-1}^{(1)}$-crystal $V^{r,s}$. Denote the virtualization by
$v:B^{r,s} \hookrightarrow V^{r,s}$. We perform explicit
computations using the tableau realization of $U_q(A_{2n-1})$-crystals
in~\cite{KN:1994} and $0$-arrows given by~\cite{S:2002}. We have
\begin{equation*}
\begin{split}
v(b)&=(2n-k)^s \cdots (r+2)^s (r+1)^s k^s \cdots 2^s 1^s \otimes r^s
\cdots 2^s 1^s \\
v(y)&=(2n)^s (2n-k-1)^s \cdots (r+2)^s (r+1)^s (k+1)^s \cdots 3^s
2^s \otimes
r^s \cdots 2^s 1^s\\
v(f_0^s y)&=(2n-k-1)^s \cdots (r+2)^s (r+1)^s (k+1)^s \cdots 2^s 1^s
\otimes r^s \cdots 2^s 1^s \\
&= v(u_{s\om_{k+1}}).
\end{split}
\end{equation*}
\end{proof}

The next theorem is the main result of this paper. It shows that under the
isomorphism between the Demazure and the KR crystals as classical
crystals zero arrows map to zero arrows. In addition it yields the
isomorphism~\eqref{E:affiso} without the assumption that the KR crystal
$B^{r,c_r s}$ is perfect.

\begin{theorem}\label{thm:Dem KR} Let $(r,s)\in \Ih\times \Z_{>0}$.
Suppose that $r\in\IS$, or $c_r\om_r=\theta$, or $\geh$ is of
nonexceptional affine type. Write $t_{-c_r\om_r^*} = w \tau$ with
$w\in W$ and $\tau\in\Sigma$. Then there is an affine crystal
isomorphism
\begin{equation} \label{E:affiso}
\begin{split}
  B(s \La_{\tau(0)}) &\cong B^{r,c_r s} \otimes B(s\La_0) \\
  u_{s\La_{\tau(0)}} &\mapsto u' := u \otimes u_{s\La_0}
\end{split}
\end{equation}
where $u$ is the element specified by Assumption~\ref{A:KR}~(\ref{A:u}). It
restricts to an isomorphism
\begin{equation} \label{E:Dem}
  D(c_r\om_r,s) \cong B^{r,c_r s} \otimes u_{s\La_0}
\end{equation}
where both sides of \eqref{E:Dem} are regarded as full subcrystals
of their respective sides in \eqref{E:affiso}.
\end{theorem}
\begin{proof}
Let $w_2$ be the minimum length coset representative in $\Wc w$.
Then $w=w_1w_2$ is a length-additive factorization with
$w_1=ww_2^{-1}\in \Wc$. We choose a reduced word of $w$ by
concatenating reduced words of $w_1$ and $w_2$. We claim that it
suffices to establish the following assertions.
\begin{enumerate}
\item[(A1)] There is a bijection
\begin{equation} \label{E:B'}
\begin{split}
  B_{w_2}(s\La_{\tau(0)}) &\rightarrow B' := f_{w_2}(u') \\
  u_{s\La_{\tau(0)}} &\mapsto u'
\end{split}
\end{equation}
that preserves all arrows in $f_{w_2}$.
\item[(A2)] $B'\subset B^{r,c_r s} \otimes u_{s\La_0}$.
\end{enumerate}
Suppose (A1) and (A2) hold. Since $w_1\in \Wc$,
$B_{w_2}(s\La_{\tau(0)})$ contains all the classical highest weight
vectors of $D(c_r\om_r,s)$. By (A1) these classical highest weight
vectors correspond to the classical highest weight vectors in $B'$.
Let $B''\subset B^{r,c_rs} \otimes B(s\La_0)$ be the classical
subcrystal generated by $B'$; by (A2) $B'' \subset B^{r,c_rs}
\otimes u_{s\La_0}$. By Demazure theory for highest weight modules
over simple Lie algebras, the bijection \eqref{E:B'} extends
uniquely to a classical crystal isomorphism $D(c_r\om_r,s) \cong
B''$. By Assumption \ref{A:KR} and Remark \ref{R:KRDem}
we have $B'' = B^{r,c_rs} \otimes u_{s\La_0}$. So we have a bijection
\begin{align} \label{E:Dembij}
  D(c_r\om_r,s) \cong B^{r,c_rs} \otimes u_{s\La_0}
\end{align}
which is an isomorphism of classical crystals that extends the
bijection \eqref{E:B'}. It follows that $B^{r,c_rs} \otimes
u_{s\La_0}$ and therefore $B^{r,c_rs} \otimes B(s\La_0)$, have a
unique affine highest weight vector, namely, $u'$. By \cite[Prop.
2.4.4]{KKMMNN:1992} there is an affine crystal isomorphism
\eqref{E:affiso}. It must extend the bijection \eqref{E:Dembij}, and
the Theorem follows.

We prove (A1) and (A2) by cases.

If $r\in \IS$ then by \eqref{E:specialWSigma} $w_2$ is the identity,
$B_{w_2}(s\La_{\tau(0)})=\{u_{s\La_{\tau(0)}}\}$, $B'=\{u'\}$,
$c_r=1$, and $B^{r,s} \cong B(s\om_r)$ as a classical crystal with
classical highest weight vector $u$. In this case (A1) and (A2) are
immediate. This is the only case where $\om_r^*\ne \om_r$.

If $c_r\om_r=\theta$
then $\tau$ is the identity, $w_1=s_\theta$ and $w_2=s_0$. By
Assumption~\ref{A:KR}~(\ref{A:u}), $B_{w_2}(s\La_0)$ and $B'$ are the $0$-strings
of $u_{s\La_0}$ and $u'$ respectively. The elements are at the
dominant ends of their respective $0$-strings, which both have
length $s$. This gives (A1). (A2) follows by the signature rule and
Assumption~\ref{A:KR}~(\ref{A:u}).

Otherwise we assume that $\geh$ is of nonexceptional affine type and
$r\in I\setminus\IS$. Then $w_2$ is given in Proposition
\ref{P:trans}. We use the notation of Remark \ref{R:w1w2} throughout
the rest of the proof. Since $K\subsetneq I$, $\geh_K$ is a simple
Lie algebra and Assumption~\ref{A:KR}~(\ref{A:regular}) implies that $B^{r,c_rs}$
decomposes into a direct sum of $K$-components, each of which is
isomorphic to the crystal graph of an irreducible highest weight
module for $U_q(\geh_K)$. We have the $K$-crystal isomorphisms
\begin{align} \label{E:Kiso}
  B_{w_2}(s\La_{\tau(0)}) \cong B_{w_2}(s\omt_{\tau(0)}) =
  B(s\omt_{\tau(0)}) \cong B'.
\end{align}
The first isomorphism holds by restriction from an $I$-crystal to a
$K$-crystal. The equality holds by Remark \ref{R:w1w2} and Demazure
theory for the simple Lie algebra $\geh_K$. We have
$B_{w_2}(s\omt_{\tau(0)})\cong B'$, since both sides are generated
by $f_{w_2}$ (with $w_2\in W_K$) applied to $K$-highest weight
vectors of $K$-weight $s\omt_{\tau(0)}$; see Assumption~\ref{A:KR}~(\ref{A:u}).
This establishes (A1).

For types $D_n^{(1)}, B_n^{(1)}, A_{2n-1}^{(2)}$ we have $c_r=1$ for
all $r$ and $\tau=\tau_0$ or $\tau=\tau_1$ (and $\tau(0)=0$ or
$\tau(0)=1$) according as $r$ is even or odd. Here
$u=u_{s\om_{\tau(0)}}\in B(s \om_{\tau(0)})\subset B^{r,c_rs}$,
where $\om_0=0$ by convention.

We consider the decomposition of $B^{r,c_r s}$ into $K$-components,
which we call $D_r$-components. Note that $0$ and $1$ are the spinor
nodes in $D_r$. Now $u_{c_r s\om_r}\in B^{r,c_rs}$ is a $D_r$-lowest
weight vector of $D_r$-weight $-2s\omt_0$. Therefore there is a
$D_r$-crystal embedding
\begin{equation*}
\begin{split}
  B(2s\omt_0) \otimes u_{s\La_0} &\rightarrow
  B(s\omt_{\tau(0)})^{\otimes 2} \otimes B(s\La_0) \\
  u_{c_r s\om_r} \otimes u_{s\La_0} &\mapsto
  u_{-s\omt_0}^{\otimes 2} \otimes u_{s\La_0}.
\end{split}
\end{equation*}
But by Lemma \ref{lem:y} there is a $D_r$-path from $u'$ to $u_{c_r
s\om_r} \otimes u_{s\La_0}$ that never changes the right hand tensor
factor. Therefore there is a $D_r$-embedding
\begin{equation*}
\begin{split}
  B' &\rightarrow B(s\omt_{\tau(0)})^{\otimes 2} \otimes B(s\La_0)
  \\
 u'= u \otimes u_{s\La_0} &\mapsto u_{s\omt_{\tau(0)}} \otimes u_{-s\omt_0} \otimes
  u_{s\La_0}.
\end{split}
\end{equation*}
The image of $u'$ is uniquely determined by Assumption~\ref{A:KR}~(\ref{A:u})
since $u_{s\omt_{\tau(0)}} \otimes u_{-s\omt_0}$ is the unique
element of $B(s\omt_{\tau(0)})^{\otimes 2}$ with $\ve=s\La_0$ and
$\vp=s\La_{\tau(0)}$.

The form of the image of $u'$ now clearly shows that when $f_{w_2}$
is applied to $u'$ it only acts on the left hand tensor factor. This
implies (A2).

Next let us consider type $C_n^{(1)}$ for $r\notin\IS$; for such
$r$, $c_r=2$ and $\tau$ is the identity. Here $u$ is the unique
element in the one-dimensional $C_n$-crystal in $B^{r,2s}$.
We decompose $B^{r,2s}$ as a $K$-crystal, which is a $C_r$-crystal in this
case. All other arguments go through as for type $D_n^{(1)}$.

Types $D_{n+1}^{(2)}$ and $A_{2n}^{(2)}$ follow in the same fashion.
In this case the decomposition of $B^{r,c_rs}$ as a $K$-crystal is a
$B_r$ crystal.
\end{proof}

\begin{remark}\label{rem:missing_nodes}
We expect Theorem \ref{thm:Dem KR} to hold for any affine algebra
$\geh$ and any Dynkin node $r\in\Ih$.
Our proof requires a special property, that the minimum length coset
representative $w_2$ of Proposition \ref{P:trans} has a certain
form, namely, in the notation of \eqref{E:wwla}, $w_2=\ww{\la}$
where $\la$ is a fundamental weight for some subalgebra $\geh_K$
where $K\subsetneq I$. This property of $w_2$ does not hold for the
trivalent node in type $E_6^{(1)}$. For such nodes a different
strategy is required.
\end{remark}

\begin{remark} \label{R:sigma}
In the notation of Lemma \ref{L:F} we expect that for any affine
algebra $\geh$ with affine Dynkin diagram $X$ and any
$\sigma\in\Aut(X)$, there is a bijection $\sigma:B^{r,c_r s}
\rightarrow B^{\sigma'(r),c_r s}$ such that \eqref{E:sigmastring}
holds. In particular, for any $\sigma\in \Aut(X)$, we expect that
there is an automorphism $\sigma$ on $B^{r,c_r s}$ satisfying
\eqref{E:sigmastring} if and only if $\sigma'(r)=r$. By Remark
\ref{R:autcl} this means that every special Dynkin automorphism
$\sigma\in \Sigma$ should induce an automorphism of each
$B^{r,c_rs}$. In contrast, for the nonspecial automorphism
$\sigma=(0,1)$ of $D_n^{(1)}$, $\sigma'=(n-1,n)$ is not the identity
and $\sigma$ induces a bijection $B^{n-1,s}\rightarrow B^{n,s}$
satisfying \eqref{E:sigmastring}.
\end{remark}
Remark~\ref{R:sigma} comes into play in Section~\ref{S:connect} and the following Theorem.

\begin{theorem} \label{thm:sigma} For the cases in
Assumption~\ref{A:KR}~(\ref{A:auto}) where $\sigma$ is defined, there exist
unique maps
\begin{equation*}
\Psi: D(\omega_r,s) \hookrightarrow  B^{r,c_rs} \text{ and }
\Psi^\sigma: D^\sigma(\omega_r,s) \hookrightarrow B^{r,c_rs}.
\end{equation*}
The maps are induced by $\Psi(u_{s \La_0}) = u$ and
$\Psi^\sigma(u_{s\La_{\sigma(0)}}) = \sigma(u)$.
\end{theorem}
\begin{proof} The map $\Psi^\sigma$ is obtained by applying $\sigma$
to everything in sight.
\end{proof}

\begin{corollary} \label{cor:KR unique}
The affine structure of $B^{r,c_rs}$ is uniquely determined.
\end{corollary}

\begin{theorem} \label{T:KRtensor} Suppose that $\la=\sum_{r\in\Ih}
m_r c_r \om_r$ with $m_r\in \Z_{\ge0}$ and $m_r>0$ only when $r$ is
as in Theorem~\ref{thm:Dem KR}. Write $t_{-\la^*} = w\tau$ for $w\in
W$ and $\tau\in\Sigma$. Assume that for each $k\in\IS$ and every
$r\in\Ih$ with $m_r>0$, the special Dynkin automorphism $\tau_k\in
\Sigma$ induces an automorphism of $B^{r,c_rs}$ that sends
$i$-arrows to $\tau_k(i)$-arrows. Then for every $r'\in\IS$ there is
an isomorphism
\begin{align*}
  B(s\La_{\tau(r')}) \cong (\bigotimes_{r\in\Ih} (B^{r,c_rs})^{\otimes
  m_r})
  \otimes B(s\La_{r'})
\end{align*}
which restricts to an isomorphism of full subcrystals
\begin{align*}
  B_{\tau_{r'}^{-1}w\tau_{r'}}(s\La_{\tau(r')}) \cong (\bigotimes_{r\in\Ih} (B^{r,c_rs})^{\otimes m_r})
  \otimes u_{s\La_{r'}}.
\end{align*}
\end{theorem}
\begin{proof} Induction allows a straightforward reduction to the
case of one KR tensor factor. Applying a special Dynkin automorphism
allows the reduction to the case $r'=0$, which is Theorem
\ref{thm:Dem KR}.
\end{proof}

\begin{corollary}
Let $\la$ be as in Theorem~\ref{T:KRtensor}. Then the Demazure
crystal $D(\la,s)$ can be extended to a full affine crystal by
adding 0-arrows.
\end{corollary}

\begin{remark}
This proves Conjecture 1 in~\cite{FL:2004} on the level of crystals.
However it is not yet clear whether there exists a global basis of
the Demazure module, whose corresponding crystal basis is the one
given in Theorem~\ref{T:KRtensor}. For level $s=1$, 
Theorem~\ref{T:KRtensor} was proved using the Littelmann path
model in~\cite[Proposition 3]{FL:2005}.
\end{remark}

\section{Reaching the classical highest weight vectors of a KR
crystal} \label{S:hw}

In the proof of Lemma \ref{lem:y}, explicit paths in the KR crystal
were given, from the element $u$ to certain classical highest weight
vectors in the KR crystal. For $\geh$ of nonexceptional affine type
and for each KR crystal $B^{r,c_r s}$, we shall give (without proof)
an explicit way to reach each classical highest weight vector in
$B^{r,c_rs}$ from the element $u$ of Assumption \ref{A:KR}.

If $r\in \IS$ then the KR crystal $B^{r,c_r s}$ is connected as a
classical crystal and the problem is trivial. This includes all
$r\in\Ih$ for $A_n^{(1)}$.

So we now assume $r\not\in\IS$.

We shall use the standard realizations of the weight lattices of
$B_n,C_n,D_n$ by sublattices of $((1/2)\Z)^n$. We let
$\om_i=(1^i,0^{n-i})$ for $i\in\Ih$ nonspin. Since $r\not\in\IS$ the
only spin weight we need is $\om_n=(1/2)(1^n)$ in type $B_n$, and in
that case $c_n=2$. Thus all the weights we must consider, correspond
to partitions, elements in $\Z_{\ge0}^n$ consisting of weakly
decreasing sequences. Moreover, for the nonexceptional affine
algebras the KR crystals are multiplicity-free as classical
crystals.

For $\geh$ of type $B_n^{(1)}$, $D_n^{(1)}$, or $A_{2n-1}^{(2)}$,
$B(\la)$ occurs in $B^{r,c_rs}$ if and only if the diagram of the
partition corresponding to $\la$, is obtained from the $r\times s$
rectangular partition by removing vertical dominoes. Let $t=0$ or
$t=1$ according as $r$ is even or odd. We have
\begin{align*}
u_\la = \left( \prod_{i=(r-t)/2}^1 f_0^{\la_{2i}} (f_2^{\la_{2i}}
f_3^{\la_{2i}} \dotsm f_{2i-1+t}^{\la_{2i}})
(f_1^{\la_{2i}}f_2^{\la_{2i}} \dotsm f_{2i-2+t}^{\la_{2i}}) \right)
u
\end{align*}
where the product is formed from left to right using decreasing
indices $i$.
\begin{example}
Let $\geh$ be of type $D_7^{(1)}$, $(r,s)=(5,4)$ and $\la$ be the
weight $\om_5+\om_3+2\om_1$. Then $t=1$, $\la$ is the partition
$(4,2,2,1,1)$, and the sequence of lowering operators is
$(f_0f_2f_3f_4f_1f_2f_3)(f_0^2f_2^2f_1^2)$. This is applied to the
classical highest weight vector of weight given by the partition
$(4)$, and the parenthesized subexpressions successively yield
classical highest weight vectors corresponding to the partitions
$(4,2,2)$, and $(4,2,2,1,1)$ respectively.
\end{example}

For $\geh$ of type  $C_n^{(1)}$, $A_{2n}^{(2)}$ or $D_{n+1}^{(2)}$,
the partitions corresponding to classical highest weights in
$B^{r,c_rs}$ are precisely those of the form $c_r \la =
(c_r\la_1,c_r\la_2,\dotsc)$ where $\la$ runs over the partitions
contained in the $r\times s$ rectangle. We have
\begin{align*}
u_{c_r\la} = \left(\prod_{i=r}^1 f_0^{c_r\la_i} f_1^{c_r \la_i}
\dotsm f_{i-1}^{c_r\la_i}\right) u
\end{align*}
where the product of operators is formed from left to right as $i$
decreases.
\begin{example}
Let $\geh$ be of type $C_3^{(1)}$, $(r,s)=(2,3)$, and
$\la=\om_2+2\om_1$. Then we have $c_r=2$, the partition $\la=(3,1)$,
and the sequence of lowering operators $(f_0^2f_1^2)(f_0^6)$. This
is applied to the classical highest weight vector of weight $0$
(corresponding to the empty partition). After $f_0^6$ the classical
weight is given by the partition $(6)$ and after $f_0^2f_1^2$ one
has the partition $(6,2)=2\la$.
\end{example}

\section{Connectedness} \label{S:connect}

Theorem~\ref{thm:Dem KR} shows that the KR crystals $B^{r,c_rs}$ are connected.
In this section we show that the tensor product of two KR crystals is also connected
by providing an algorithm which for any given element in the crystal yields a string
of operators $e_i$ (or $f_i$) to reach a given special element. This algorithm is also
useful in defining crystal morphisms such as the combinatorial $R$-matrix.
Since KR crystals and their tensor products are not highest weight crystals,
it is not completely obvious which sequence of raising operators $e_i$
will yield a given special element.

Here we give a construction on how to
reach $u_1 \otimes u_2 \in B^{r_1,c_{r_1}s_1} \otimes B^{r_2,c_{r_2}s_2}$ where $u_1$ is
the unique elements of $B^{r_1,c_{r_1}s_1}$ with $\ve(u_1)=s_1 \La_0$ and
$\vp(u_1)=s_1 \La_{\tau_1(0)}$ as required in Assumption~\ref{A:KR}~(\ref{A:u}),
and $u_2$ is the unique element in $B^{r_2,c_{r_2}s_2}$ with
$\ve(u_2)=s_2 \La_{\tau_2^{-1}(0)}$ and $\vp(u_2)=s_2\La_0$ as required
in Assumption~\ref{A:KR}~(\ref{A:u}) and Remark~\ref{R:sigma}.

By Theorems~\ref{thm:Dem KR} and~\ref{thm:sigma} we have the following isomorphism
of affine crystals
\begin{equation*}
\begin{split}
B^{r_1,c_{r_1}s_1} \otimes B^{r_2,c_{r_2}s_2} \otimes B(s_2 \Lambda_{\tau^{-1}_2(0)})
&\cong B^{r_1,c_{r_1}s_1} \otimes B(s_2 \Lambda_0)\\
u_1 \otimes u_2 \otimes u_{s_2 \La_{\tau_2^{-1}(0)}} &\mapsto
u_1 \otimes u_{s_2 \La_0}.
\end{split}
\end{equation*}
Assume that $s_1\ge s_2$. Acting with raising operators $e_i$ with
$i\in I$ one can bring any element $b_1 \otimes b_2 \otimes u_{s_2
\Lambda_{\tau_2^{-1}(0)}}$ into the form $c_1 \otimes u_2 \otimes
u_{s_2 \Lambda_{\tau_2^{-1}(0)}}$ since by the tensor product rule
the $e_i$ will eventually act on the right tensor factors and by
Theorem~\ref{thm:Dem KR} $b_2\otimes u_{s_2
\Lambda_{\tau_2^{-1}(0)}}$ is connected to $u_2 \otimes u_{s_2
\Lambda_{\tau_2^{-1}(0)}}$ . Once such an element is reached, tensor
from the right by $u_{(s_1-s_2)\La_0}\in B((s_1-s_2)\La_0)$ to
obtain
\begin{multline*}
B^{r_1,c_{r_1}s_1} \otimes B^{r_2,c_{r_2}s_2} \otimes B(s_2 \Lambda_{\tau_2^{-1}(0)})
\otimes B((s_1-s_2)\Lambda_0)\\
\cong B^{r_1,c_{r_1}s_1} \otimes B(s_2 \Lambda_0) \otimes
B((s_1-s_2)\Lambda_0)
\end{multline*}
under which $c_1 \otimes u_2 \otimes u_{s_2
\Lambda_{\tau_2^{-1}(0)}} \otimes u_{(s_1-s_2)\La_0}$ maps to $c_1
\otimes u_{s_2\La_0} \otimes u_{(s_1-s_2)\La_0}$. The latter element
is the image of the vector $c_1 \otimes u_{s_1\La_0}$ under the
embedding of affine crystals $B^{r_1,c_{r_1}s_1} \otimes B(s_1\La_0)
\rightarrow B^{r_1,c_{r_1}s_1} \otimes B((s_1-s_2)\La_0) \otimes
B(s_2\La_0)$.

Now from $c_1\otimes u_{s_1\Lambda_0} \in B^{r_1,c_{r_1}s_1} \otimes B(s_1
\Lambda_0)$ one can reach $u_1\otimes u_{s_1\Lambda_0}$ using $e_i$ with $i\in I$.

If $s_1<s_2$ we tensor from the left with the dual crystals.
Explicitly,
\begin{equation*}
B^\vee(s_1 \Lambda_{\tau_1(0)}) \otimes B^{r_1,c_{r_1}s_1} \otimes
B^{r_2,c_{r_2}s_2} \cong B^\vee(s_1 \Lambda_0) \otimes B^{r_2,c_{r_2}s_2}.
\end{equation*}
The lowest weight element $u^\vee_{s_1 \Lambda_0}\in B^\vee(s_1
\Lambda_0)$ corresponds to $u^\vee_{s_1 \Lambda_{\tau_1(0)}} \otimes
u_1 \in B^\vee(s_1 \Lambda_{\tau_1(0)}) \otimes B^{r_1,c_{r_1}s_1}$.
Acting with lowering operators $f_i$ with $i\in I$ one can bring any
element $u^\vee_{s_1\Lambda_{\tau_1(0)}} \otimes b_1 \otimes b_2$
into the form $u^\vee_{s_1\Lambda_{\tau_1(0)}} \otimes u_1 \otimes
c_2$. Once this element is reached, tensor on the left by
$u^\vee_{(s_2-s_1)\La_0}\in B^\vee((s_2-s_1)\Lambda_0)$, obtaining
the element $u^\vee_{(s_2-s_1)\La_0} \otimes
u^\vee_{s_1\Lambda_{\tau_1(0)}} \otimes u_1 \otimes c_2$, which can
be identified with $u^\vee_{s_2\La_0} \otimes c_2 \in
B^\vee(s_2\La_0) \otimes B^{r_2,c_{r_2}s_2}$. Now move down to the
lowest weight vector $u^\vee_{s_2 \Lambda_0} \otimes u_2$ using
$f_i$ with $i\in I$.

As a result of the above construction we obtain the following corollary:
\begin{corollary} \label{cor:connect}
The tensor product $B^{r_1,c_{r_1}s_1} \otimes B^{r_2,c_{r_2}s_2}$ of KR crystals
is connected.
\end{corollary}

The combinatorial $R$-matrix is a crystal morphism. More precisely
\begin{equation*}
 R:B^{r_1,c_{r_1} s_1} \otimes B^{r_2,c_{r_2} s_2}
  \to B^{r_2,c_{r_2} s_2} \otimes B^{r_1,c_{r_1} s_1}
\end{equation*}
satisfies $R\circ e_i = e_i \circ R$ and $R\circ f_i = f_i \circ R$ for all $i\in I$.
There exists a unique element $u_{c_{r_k} s_k \om_{r_k} }\in B^{r_k,c_{r_k} s_k}$
and by weight considerations $R$ must map
$R(u_{c_{r_1}s_1 \om_{r_1}} \otimes  u_{c_{r_2}s_2 \om_{r_2}})=
u_{c_{r_2}s_2 \om_{r_2}} \otimes  u_{c_{r_1}s_1 \om_{r_1}}$.
Assume  that $s_1\ge s_2$. Then for any element $b_1  \otimes b_2 \in B^{r_1,c_{r_1} s_1} \otimes
B^{r_2,c_{r_2} r_2}$ the above algorithm provides a sequence
$e_{\{i\}} := e_{i_1} e_{i_2} \cdots e_{i_\ell}$ such that $e_{\{i\}}(b_1\otimes b_2)= u_1\otimes u_2$.
In particular, $e_{\{j\}}(u_{c_{r_1}s_1 \om_{r_1}} \otimes  u_{c_{r_2}s_2 \om_{r_2}})=
u_1\otimes u_2$.  Set $f_{\{\leftarrow i\}} := f_{i_\ell} \cdots f_{i_1}$. Then
\begin{equation*}
R(b_1 \otimes b_2)= f_{\{\leftarrow i\}} e_{\{j\}}
 (u_{c_{r_2}s_2 \om_{r_2}} \otimes  u_{c_{r_1}s_1 \om_{r_1}}).
\end{equation*}
For the case $s_1<s_2$ a similar construction works where $f_i$ and $e_i$ are
interchanged.

\end{document}